# Quasiperiodic music

Darren C. Ong

Department of Mathematics, Xiamen University Malaysia, Sepang, Malaysia

**ABSTRACT**

Using the definition of quasiperiodic function as a motivation, we introduce the idea of quasiperiodic music and detail the composition process of a quasiperiodic music piece, *Raindrops in A minor*. We also discuss connections between quasiperiodic music and other works of music theory and composition that make use of aperiodic order or periodic order with large periods, such as Lindenmayer systems, Vuza canons, Messiaen's *Quatuor pour la Fin du Temps*, and the phase music of Steve Reich.

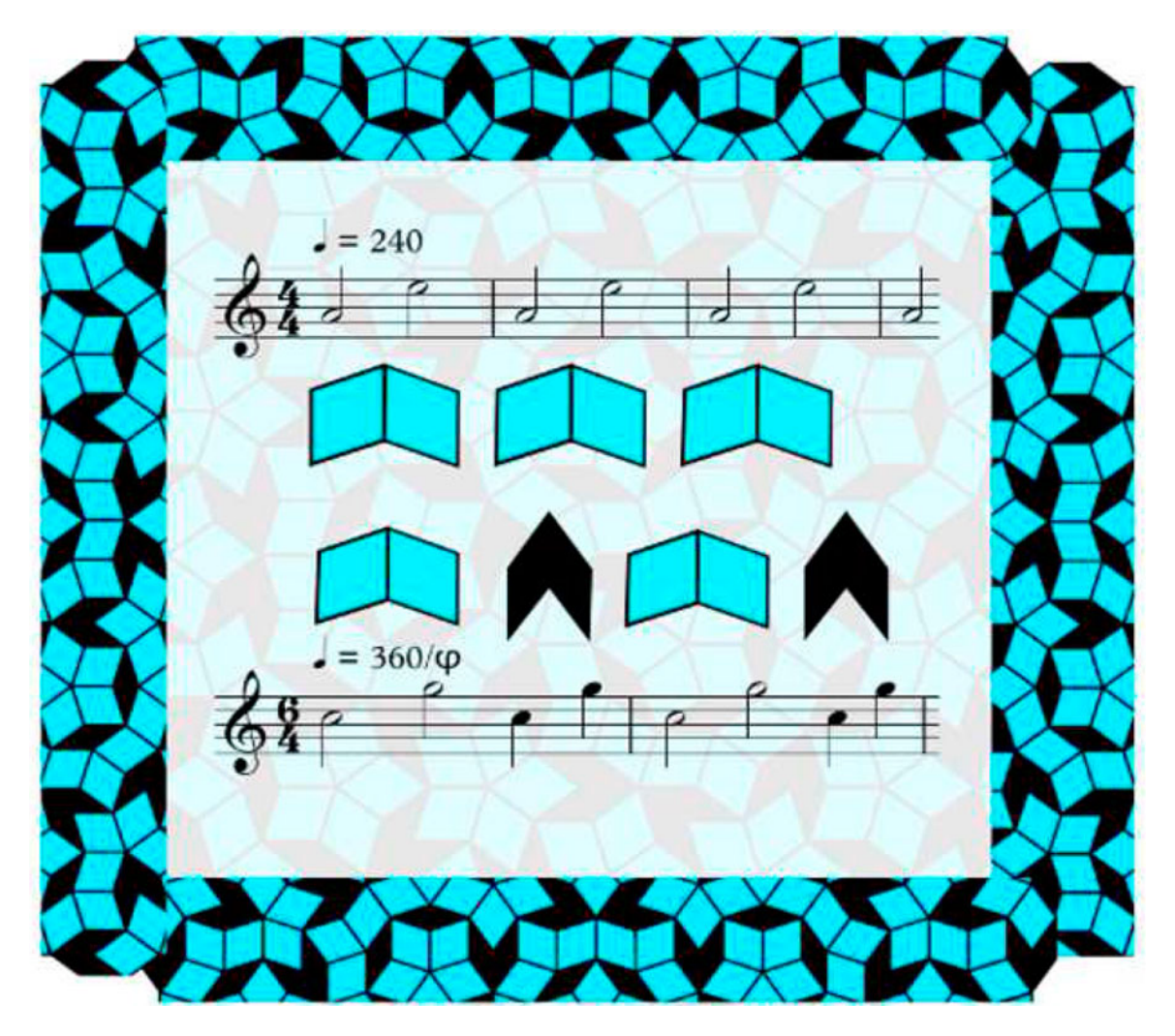

## 1. Introduction

In recent decades, there has been a surge of interest in aperiodic order. 'Aperiodic order' is a rather broad term that encompasses both discrete and continuous structures. These different structures are loosely grouped under one umbrella term because they

**CONTACT** Darren C. Ong darrenong@xmu.edu.my

are aperiodic, at the same time, each of them is close to periodic in some mathematically precise way (although the sense in which they are 'close' may differ). See Baake and Grimm (2013, 2017), Baake et al. (2016), Kellendonk et al. (2015), Barber (2008) and Dal Negro (2013) for recent books detailing the impact of aperiodic structures in pure mathematics, mathematical physics, condensed matter physics, and optics.

As an example of a discrete aperiodic pattern, we take the Fibonacci word, which is an infinite binary string starting with 0100101001001. This is listed as sequence A003849 in the Online Encyclopedia of Integer Sequences (Sloane, 2019). Note that the Fibonacci *word* is different from the more well-known Fibonacci *sequence* $1, 1, 2, 3, 5, 8, 13, \ldots$, in which every term in the sequence (after the second term) is the sum of the previous two. To form the Fibonacci word, we start with initial bit strings $S_0 := 0$ and $S_1 := 01$, and then we form a sequence of bit strings for which every string in the sequence is obtained by concatenating the preceding two terms. That is, $S_{n+2} := S_{n+1} \odot S_n$ for $n = 0, 1, 2, \ldots$. The first few strings in the sequence will be:

$$\begin{aligned} S_0 &= 0 \\ S_1 &= 01 \\ S_2 &= 010 \\ S_3 &= 01001 \\ S_4 &= 01001010 \end{aligned}$$

and it is clear that the strings in the sequence will converge to the infinite string 0100101001001 . . ., which we call the Fibonacci word.

This Fibonacci word contains exactly $n+1$ different $n$-bit substrings. We define a complexity function on a binary string as follows. Let $\sigma(n)$ be the number of substrings of length $n$ of that string. We can observe that, for the Fibonacci word, $\sigma(n) = n + 1$. For a periodic binary string, it is easy to see that $\sigma(n)$ is bounded. By the Morse-Hedlund Theorem (Baake & Grimm, 2013, Proposition 4.11), for any aperiodic infinite string $\sigma(n) \geq n + 1$. It is in this sense that we can say that the Fibonacci word is aperiodic but very close to periodic.

Let us now consider an example in the continuous setting, with almost periodic functions. These are explained in detail in Appendix 1 of Avron and Simon (1981), that we briefly summarize here. Given a function $f$, let us define the translated function $f_t(x) = f(x - t)$. A (*Bochner*) *almost-periodic function* is a bounded continuous function on $\mathbb{R}^\nu$ such that the set of functions $\{f_t | t \in \mathbb{R}^\nu\}$ form a precompact set with respect to the supremum norm (a precompact set is a set whose closure is compact). The function $f$ is periodic if and only if $\{f_t | t \in \mathbb{R}^\nu\}$ is a compact set itself (Remling, 2019), so we see again that an almost-periodic function is close to a periodic one.

Fibonacci words and almost periodic functions are two of the most prominent examples of aperiodic order, but there are many others. Perhaps, the most visually iconic example of aperiodic order would be the Penrose tiling (Figure 1), which is an aperiodic tiling of the plane. We say that a tiling is periodic if it is invariant under a translation. Imagine an infinite floor tiled with $1 \times 1$ inch square tiles. If we were to shift every tile 1 inch to the left, the pattern of the tiling would be unchanged. In this sense, this square tiling is periodic. The Penrose tiling, in constrast, is not periodic because no translation will leave the tiling perfectly invariant in the same way. Curiously, we can also see the Penrose tiling as close to periodic in the following way: any finite portion of the Penrose tiling appears infinitely many times. These patterns are very compelling visually, and indeed they have appeared in the visual arts and architecture dating back almost a millennium ago.

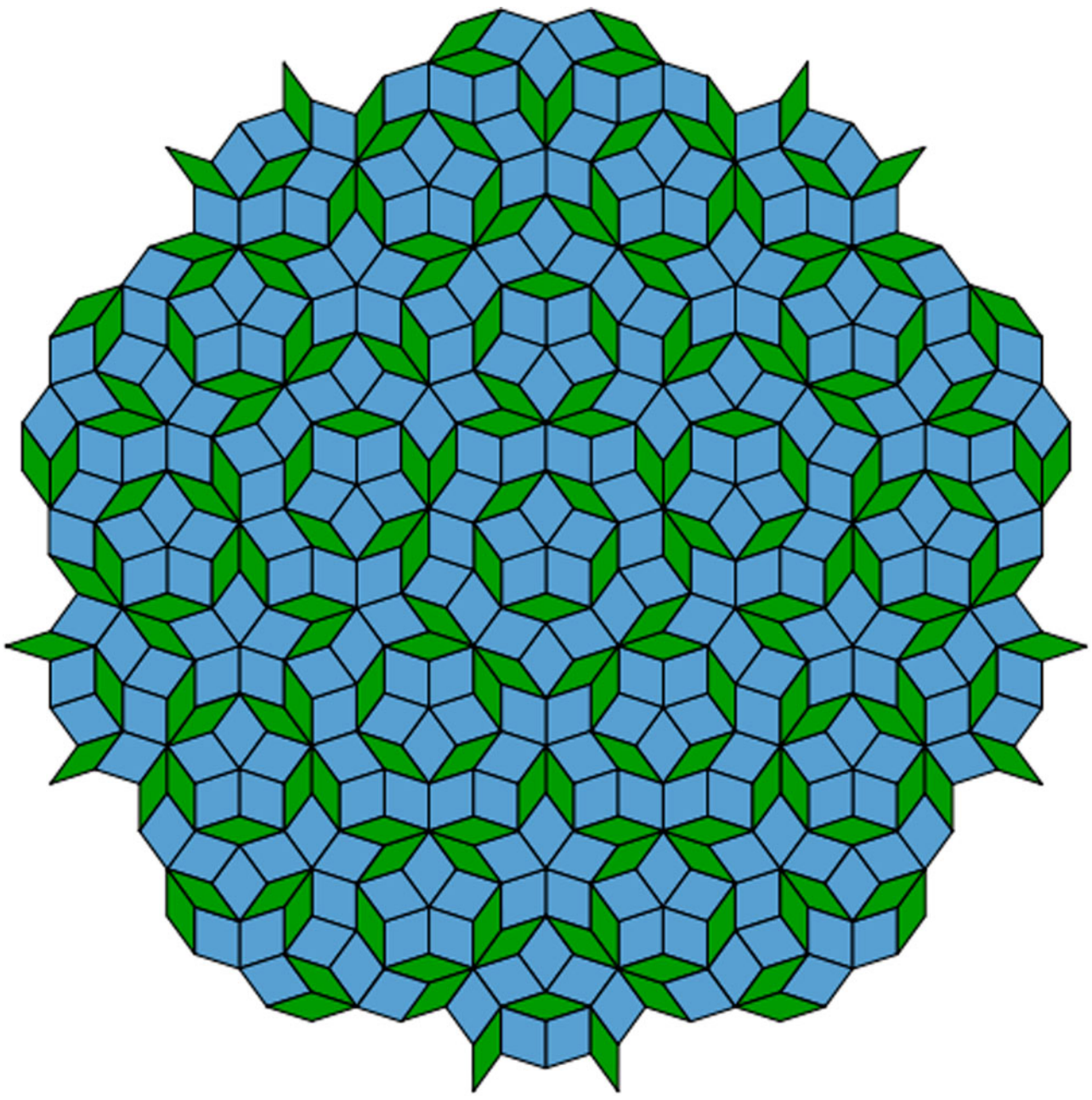

**Figure 1.** A small section of the Penrose tiling of the P3 Rhombus type. This image was created by User:Inductiveload from Wikipedia and released under the public domain.

Scientists have realized that these patterns are also of great significance in physics and chemistry, especially after the 1984 discovery of the quasicrystal, which is a material with an aperiodically-ordered atomic structure. Since then, there have been many efforts in physics and mathematics to model quasicrystals and to understand aperiodically-ordered mathematical structures. One of the most important ways to model quasicrystalline structures exploits a *quasiperiodic function*. We will define quasiperiodic functions precisely later in Definition 3.1; for now, we can informally say that these functions are generated by combining periodic functions of different frequencies that never match up precisely. The quasiperiodic functions are a special case of almost periodic functions, as explained in Theorem A.1.2 of of Avron and Simon (1981).

Aperiodically ordered patterns have become very prominent in the visual arts, in science, and in mathematics, so it is unsurprising that they have also inspired musicians and music theorists. In particular, musical rhythm has long been understood to have a rich mathematical structure (Sethares, 2014). Discrete aperiodic order has been applied to rhythmic structures in the form of Lindenmayer systems, for instance.

In this article, we creatively explore quasiperiodic music, and we propose, as an application, the original composition *Raindrops in A minor* (Ong, 2020), which is an attempt to express musically the structure of a quasiperiodic function.

This paper is structured as follows. Section 2 discusses briefly the history of aperiodic order in the arts and sciences. Section 3 provides a formal mathematical definition of quasiperiodic functions. This is the only part of the paper that requires any knowledge of technical mathematics. A reader without a mathematical background could skip this section and still grasp the main points of the paper. In fact, the main purpose of this section is to motivate some choices made in composing *Raindrops in A minor*. Section 4 details the composition of *Raindrops in A minor*, and gives a brief technical analysis of its features. Finally, Section 5 places quasiperiodic music in the context of contemporary music scholarship. We discuss connections between our work and other studies of music theory and composition that are inspired by or related to aperiodic order.

## 2. A brief history of aperiodic order in the arts and sciences

Aperiodic order has become incredibly important in music, visual arts, physics, chemistry, and mathematics. To fully appreciate quasiperiodic music, one should maybe be familiar with the artistic and scientific trends that are underpinned by these patterns.

Aperiodically ordered structures invoke some very striking images. They have a long history in Islamic art and architecture going back to the 12th century. Examples include the Darb-i Imam shrine and the Friday Mosque in Isfahan, Iran; the walls of the courtyard of the Madrasa al-Attarin in Fez, Morocco; the Seljuk Mama Hatun Mausoleum in Tercan, Turkey; and the external walls of the Gunbad-I Kabud tomb tower in Maragha, Iran (Al Ajlouni, 2012, 2013; Lu & Steinhardt, 2007). Western interest in these patterns emerged much later. In Europe, the earliest examples of aperiodic order are probably the rhombus patterns in a 1524 doodling by Albrecht Dürer (Kemp, 2005).

In the late 20th century, the importance of quasiperiodic structures in mathematics and science became apparent. Roger Penrose discovered the tilings while investigating planar tilings. The Penrose tilings were the first aperiodic tiling of the plane using just two tiles. Dan Shechtman discovered the quasicrystal in the 1980s, an achievement that would later let him win the Nobel prize in chemistry. Following this discovery, there has been a surge of interest in mathematical models of quasicrystalline structures. As an example, Artur Avila won the Fields medal in 2014 partially thanks to his investigations on the Almost Mathieu operator, which is a Schrödinger operator whose associated potential function is quasiperiodic. The Almost Mathieu operator has an important connection to the Hofstadter butterfly, another iconic mathematical image. This connection between quasiperiodicity and the butterfly is thoroughly explored in Chapter 4 of Satija (2016).

The scientific importance of aperiodic order is not restricted to mathematical physics. Aperiodic order also appears in consdensed matter physics (Barber, 2008), optics (Dal Negro, 2013), and crystallography (Baake & Grimm, 2017).

## 3. Mathematical background

### 3.1. A mathematical definition of quasiperiodicity

There are several different ways to define a quasiperiodic function in mathematics. Here is one common definition, which is a restatement of (A.1.2) of Avron and Simon (1981).

**Definition 3.1 (Quasiperiodic function):** Let $\mathbb{S}$ be the unit circle parameterized as $[0, 2\pi)$. A function $f$ on $\mathbb{R}$ is quasiperiodic if for some positive integer $n$ there exists a continuous function $Q$ on $\mathbb{S}^n$ and a frequency vector $\omega = (\omega_1, \omega_2, \ldots, \omega_n)$ such that

$$f(x) = Q(x\omega_1, x\omega_2, \ldots, x\omega_n).$$

Technically, under Definition 3.1 all continuous periodic functions are quasiperiodic (although in the rest of this article, when we refer to a function of a piece of music as quasiperiodic, we implicitly mean that it is also aperiodic). To get an aperiodic function from the above definition, we require $n$ to be at least 2 and that the entries of the frequency vector should be rationally independent:

**Definition 3.2 (Rational Independence):** We say that a set of real numbers $\omega_1, \ldots, \omega_n$ are *rationally independent* if for rational numbers $y_i$ the only solution of the equation

$$y_1\omega_1 + y_2\omega_2 + \ldots + y_n\omega_n = 0$$

is $y_1 = y_2 = \ldots = y_n = 0$.

For example, the function

$$f(x) = \sin(x) + \sin(2\pi x) \tag{1}$$

is a quasiperiodic function that is not periodic. Using Definition 3.1, here $n = 2$, $Q(v_1, v_2) = \sin(v_1) + \sin(v_2)$ and $\omega = (1, 2\pi)$. Notice that, since $2\pi$ is an irrational number, the only rational solution to $1 \cdot y_1 + 2\pi \cdot y_2 = 0$ is $y_1 = y_2 = 0$, and thus the set $\{1, 2\pi\}$ is rationally independent.

Indeed, one simple way to obtain a quasiperiodic function is to add two periodic functions, one with a rational period and one with an irrational period. We will do something similar to create our quasiperiodic music composition in Section 4.

## 3.2. Rational approximants of irrational numbers

Even though the quasiperiodic functions we consider here are aperiodic, they can become very close to being periodic. This phenomenon is connected to the idea of *rational approximants* of irrational numbers. Definition 3.3 is given in Khinchin (1997) as a best rational approximant 'of the second kind.' We will just call it a best rational approximant because in this paper we do not consider the first kind.

**Definition 3.3 (Best rational approximant):** A fraction $a/b$ is a *best rational approximant* of a number $y$ if, for every fraction $c/d$ with $1 \leq d \leq b$ and $a/b \neq c/d$, we have $|dy - c| > |by - a|$.

Intuitively, a best rational approximant is a rational number that approximates a (usually irrational) number $y$ better than any fraction with smaller or equal denominator. For example, 22/7 is a best rational approximant of $\pi$, because it is closer to $\pi$ than any other fraction with denominator 1, 2, 3, 4, 5, 6 or 7.

**Example 3.4:** The best rational approximants of $\pi$ are 3/1, 13/4, 16/5, 19/6, 22/7, 179/57, 201/64, 223/71, 245/78, 267/85, 289/92, 311/99, 333/106, 355/113, ....

**Example 3.5:** The best rational approximants of $\varphi$, the golden ratio are fractions of consecutive Fibonacci numbers, that is, $2/1, 3/2, 5/3, 8/5, /13/8, \ldots$.

These best rational approximants of a number $y$ are derived from the continued fraction expansion of $y$. More precisely, the best rational approximants are the 'convergents' and 'semi-convergents' of the continued fraction expression of $y$. Since it is not necessary for this paper we will not explain the details; but the interested reader can check the thorough explanation in Khinchin (1997).

The best rational approximants are associated to points where a quasiperiodic function is close to being periodic. For the function in (1), $\omega_2/\omega_1 = 2\pi$. The best rational approximants of $2\pi$ then correspond to points where the two components of $f(x)$, that is, $\sin(x)$ and $\sin(2\pi x)$, are close to repeating at the same time. For example, we can use $2\pi \approx 19/3$. This leads to $(1)(19) \approx (2\pi)(3)$. In other words, 19 repetitions of $\sin(2\pi x)$ are close to 3 repetitions of $\sin(x)$. A period of $\sin(x)$ begins at $x = 0, 2\pi, 4\pi, 6\pi$ and a period of $\sin(2\pi x)$ begins at $x = 0, 1, 2, 3, \ldots, 18, 19$. The closest pair of points where the period begins from this list is $6\pi \approx 18.85$ and 19, and indeed we cannot find any pair of such points in $x \in (0, 19)$ that is closer than this pair. This is a direct consequence of the fact that 19/3 is a better approximation to $2\pi$ than any fraction with denominator 3 or less (where we define 'better' in terms of Definition 3.3). We sometimes say that 19 is a *quasiperiod* of $f(x)$.

## 4. Creating 'Raindrops in A minor'

In this section, we explain the process behind the creation of the composition *Raindrops in A minor* (Ong, 2020) and we discuss its mathematical and musical properties. This piece consists of two simple elements: a repeating measure in 4/4 time running at 240 beats per minute, and a repeating measure in 6/4 time running at $360/\varphi$ beats per minute, where $\varphi \approx 1.618$ is the golden mean. The choice of $360/\varphi$ was intended to make the tempi of the two measures approximately equal. That is, $360/\varphi \approx 222.5$ is not too far off from 240 beats per minute. We repeat both measures indefinitely, starting them at the same time. Note that $\varphi$ is an irrational number, and thus the length of two measures (in terms of time) are rationally independent. To be precise, the first measure takes one second, while the second measure takes exactly $\varphi$ seconds.

Let us now directly connect this music with a quasiperiodic function as defined in Definition 3.1. We interpret the domain $\mathbb{R}$ as time. We define the codomain to be $\mathbb{R}^4$, which represents the amplitude envelopes corresponding to the four musical notes played in this piece (that is, A4, E5, C5, G5).

In the notation of Definition 3.1, we define $Q : \mathbb{S}^2 \to \mathbb{R}^4$ in the following way. We define $Q(x\omega_1, x\omega_2) = Q_1(x\omega_1) + Q_2(x\omega_2)$, where $Q_1(x\omega_1)$ maps to the A4 and E5 components, and $Q_2(x\omega_2)$ maps to the C5 and G5 components. Note that the A4 and E5 components should be 1-periodic in $x$, since $Q_1(x\omega_1)$, the part that plays the A4 and E5 notes, repeats every second. Similarly, the C5 and G5 components repeat every $\varphi$ seconds, so $Q_2(x\omega_2)$ must be $\varphi$-periodic in $x$. However, since $Q$ is defined as a function on $\mathbb{S}^2$ parametrized as $[0, 2\pi)$, to achieve the desired periodicity we must set $\omega_1 = 2\pi$ and $\omega_2 = 2\pi/\varphi$. Given these definitions, $f(x)$ will be a quasiperiodic function that indicates the amplitudes corresponding to the four notes at a time $x$ (Figure 2 ).

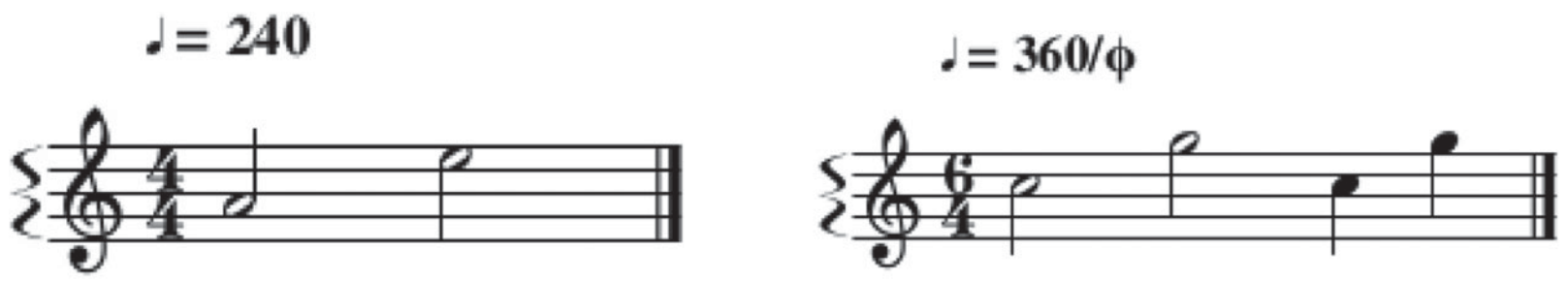

**Figure 2.** *Raindrops in A minor* is generated by playing these two infinitely repeating measures simultaneously. © Darren C. Ong.

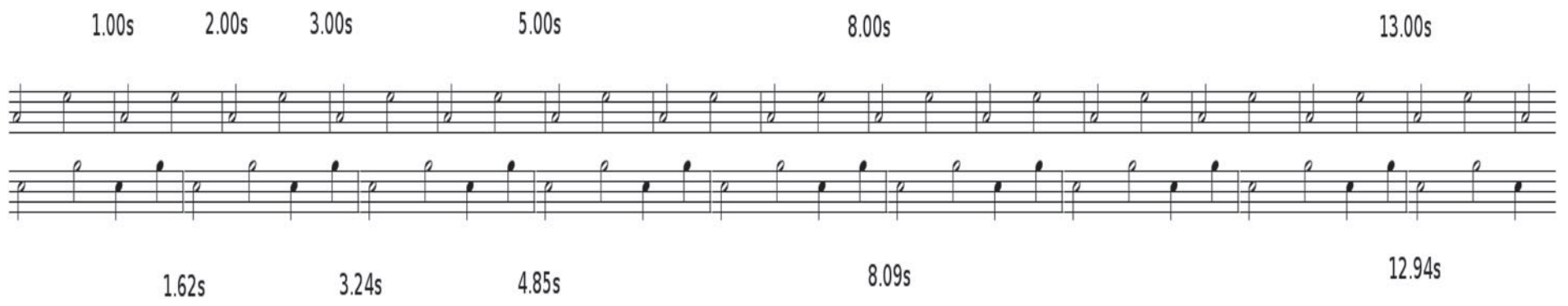

**Figure 3.** After 1,2,3,5,8 or 13 seconds, observe that the two measures in Figure 2 begin at almost exactly the same time. The larger the Fibonacci number, the closer the two measures are to starting in sync. Note that the numbers in the bottom of the figure are approximate, since $\varphi$ is an irrational number. This music is written in treble clef.

This music never repeats exactly even if it could be played for an infinite amount of time. However, there will be points in time where the music seems very close to repeating. These points are associated with rational approximants of the golden ratio. For example, since $\varphi \approx \frac{13}{8}$, this implies $13 \cdot 1 \approx 8 \cdot \varphi$. In other words, thirteen periods of the first measure takes roughly the same time as eight periods of the second measure. Indeed, we can calculate that the eighth period of the second measure ends after roughly 12.94 seconds, and the thirteenth period of the first measure obviously ends after exactly 13 seconds (see Figure 3). Following our discussion in Section 3.2, we say that this musical composition has a *quasiperiod* of 13 seconds, since the music almost repeats exactly after that amount of time. Raindrops in A minor has the pleasant property that its quasiperiods are after $1, 2, 3, 5, 8, 13, 21, \ldots$ seconds. These are precisely the Fibonacci numbers!

We can observe that the larger the quasiperiod, the closer this piece gets to repeating exactly. For example, using the approximation $\varphi \approx \frac{55}{34}$, we find that 34 periods of the second measure take about 55.01 seconds, whereas 55 periods of the first measure takes exactly 55 seconds. Using larger and larger Fibonacci numbers, we can get the two periods to synchronize as closely as we wish, but it will never synchronize exactly. This explains why this music is non-repeating, even when played for infinite time.

We should note however, that the above properties only apply for an 'idealized' version of this music. In practical terms, we created this sound file using a computer, and computers are not capable of dealing with irrational numbers like $\varphi$, resorting instead to a floating-point approximation. Thus, the full sound file we created will repeat, albeit after a large amount of time. Even if we could overcome the issues with the floating-point approximation, listeners would not be able to perceive music notes with infinite precision. In a real-world situation, a human listener would not be able to perceive that even a perfectly played version of this music is not periodic.

In terms of the music, we chose simple melodies from an A minor 7th chord. Bells seemed like a good choice of instrument, since the ringing of bells often produces notes that are played almost simultaneously. Also, a 'canon-like' musical piece played with percussion can remind us of gamelan music. In terms of software, we used Musescore to write the music, and Audacity to perform the sound mixing.

*Raindrops in A minor* was intended to be a simple example of this quasiperiodic music that demonstrates its mathematical properties. Of course, it is possible to consider more complex extensions of this idea. One could use irrational numbers other than the golden ratio $\varphi$. Different irrational numbers have different best rational approximants. This can result in a different set of quasiperiods, and so different features in the quasiperiodic music.

We could also consider writing quasiperiodic music with more parts. To create a quasiperiodic musical piece with $n$ parts for a postive integer $n$, we first find a set of $n$ rationally indepent numbers $\omega_1, \dots, \omega_n$, as in Definition 3.2. We then create music where the first part repeats every $\omega_1$ second, the second part repeats every $\omega_2$ seconds, and so on, so the $n$th part repeats every $\omega_n$ seconds. One could even imagine an entire orchestra where each instrument is playing a part with a different rationally independent period!

## 5. Related musical work

There are quite a few interesting connections between quasiperiodic music and ideas present in modern musical composition and theory. We will briefly discuss some examples in this section in order to place quasiperiodic music in the context of music scholarship.

### *5.1. Lindenmayer systems*

Lindenmayer systems (or L-systems) originated in biology as a way to model plant growth, and then inspired also musical composition. As an example, we consider an alphabet with two symbols, $A$ and $B$, with the rules $A \to AB$ and $B \to A$. Let us start with a string consisting of a single letter $A$, and then let us apply the rules to each letter in the string repeatedly, to get progressively longer strings. The first few iterations will then be $A$, $AB$, $ABA$, $ABAAB$, $ABAABABA$. The reader might notice that this is exactly the Fibonacci word we mentioned in the introduction. Indeed, these Lindenmayer systems are associated to substitution sequences, which are one of the most prominent discrete examples of aperiodic order (see Chapter 4 of Baake & Grimm, 2013).

These patterns have caught the interest of musicians. Possibily the first example of Lindenmayer systems being used in music is described in Prusinkiewicz (1986), with an algorithmic musical composition process. We quote from the conclusion of that paper:

> This paper presents a technique for generating musical scores, which consists of three steps:
> (1) A string of symbols is generated by an L-system,
> (2) This string is interpreted graphically as a sequence of commands controlling a turtle,
> (3) The resulting figure is interpreted musically, as a sequence of notes with pitch and duration determined by the coordinates of the figure segments.

> The scores generated using the above method are quite interesting. They are relatively complex (in spite of the simplicity of the underlying productions) but they also have a legible internal structure (they do not make the impression of sounds accidentally put together).

These ideas have been built upon by composers such as Hanspeter Kyburz, Enno Poppe, Phillipe Manoury, and Alberto Posada: Besada (2019) contains an interesting discussion of the impact of these Lindenmayer systems on contemporary music.

These compositions based on Lindenmayer systems are perhaps the closest antecedents of our work. Substitution sequences (which are equivalent to these Lindenmayer systems) are discrete aperiodically ordered structures, whereas quasiperiodic functions are continuous aperiodically ordered structures. The music we introduce thus shares a lot of features with the music based on Lindenayer systems. In particular, I think of Prusinkiewicz's comments on his technique producing music that is 'relatively complex in spite of simple underlying productions' and with 'legible internal structure.' This may also easily apply to quasiperiodic music. In some way, quasiperiodic music can be considered as a continuous analogue of the music based on Lindenmayer systems.

## 5.2. *Quatuor pour la Fin du Temps*

A very well-studied example of using period lengths to evoke an illusion of a music piece extending to eternity comes from Olivier Messaien's *Liturgie de Cristal*, which is the first movement of *Quatuor pour la Fin du Temps* (that is, *Quartet for the End of Time*). This piece is for violin, clarinet, cello, and piano. This movement's use of prime numbers has been well-noted in the literature. For example, on page 25 of Dingle (2016):

> 'Liturgie de Cristal' is undoubtedly the most famous example of a compositional machine in Messiaen's music. As has been discussed at length the piano part has a cycle of 29 chords combined with a set of 17 durations. Both being prime numbers, it would take 17 times 29 cycles to return to the same starting point. To this should be added the cello, which plays a cycle of 5 pitches and a sequence of 15 durations. For this machine to return to its starting point would take 4 hours 40 minutes, but Messiaen simply presents a short fragment. The effect is of a door being opened to reveal the ongoing workings, and then being closed again.

There are a lot of other works noting Messiaen's use of mathematics here, for example (Pople, 1998).

*Raindrops in A minor* differs from Messiaen's composition because it does not play with different lengths of periods in order to create an illusion of eternity; in some sense, in its idealized form, the performance of *Raindrops* is truly 'eternal.'

## 5.3. *Vuza canons*

Vuza canons are another way that aperiodicity appears in music theory. These rhythmic structures arise from combining aperiodic canons to create another canon having a large period. These canons are closely related to the problem of tiling the integers with aperiodic tiles (as in the title of Jedrzejewski, 2009). Again, there is an interplay between periodicity and aperiodicity, and complex patterns arising from repetitions of simpler ones. There is a complicated mathematical structure in these canons, including connections with group theory and cyclotomic polynomials. It is unsurprising that these canons have attracted much interest from both music theorists and mathematicians, enough so that an entire special issue of the *Journal of Mathematics and Music* is devoted to them (Andreatta & Agon, 2009). See also Andreatta (2011) for a historical survey.

### 5.4. Phasing

Phasing is a musical technique in minimal music first developed by Steve Reich, and again there are inescapable thematic connections between phase music and quasiperiodic music.

Phasing is associated with minimal music. Similar to our quasiperiodic music composition, a typical phase music composition is formed from two repetitive phrases. However, typically these two phrases are identical, and initially they are played in unison. Gradually, the phrases drift out of sync, with the second phrase lagging behind the first, intially creating an 'echo' effect, and eventually a complex combination of sounds. The phrases then drift back into unison. Note that the two phrases are typically played at approximately the same *tempo*. The slowdown to make one phrase lag behind the other is very gradual.

Probably the most famous piece of phase music is Steve Reich's *Piano Phase*. The following description of this music is in page 386 of Schwarz (1980):

> In this work two performers begin in unison playing the identical rhythmic/melodic pattern. As the first performer's pattern remains unvarying, the second pianist increases his tempo very slightly (this gradual phasing process is indicated in Reich's scores by dotted lines between measures) until he is finally one sixteenth note ahead of the unchanged figure of the first pianist. The phasing process pauses at this point, as the newly shifted rhythmic configuration is repeated several times. Soon however, the second pianist again moves slowly forward of the first, finally ending two sixteenth notes ahead of the original pattern. This sequence of gradual phase shift and repetition is repeated until the two pianists are back in unison; at this juncture the pattern changes and the whole process begins anew.

The structure of Steve Reich's phase music has received considerable mathematical analysis. In Colannino et al. (2009) and Cohn (1992) the authors perform a systematic mathematical analysis on the rhythmic structure of *Piano Phase* and its sister piece, *Violin Phase*. While these pieces show an interesting mathematical depth, they are nevertheless finite.

Conversely, quasiperiodic music is infinite and aperiodic. It never perfectly repeats, even if the pieces could be played forever. Despite these differences, it is easy to see thematic similarities between quasiperiodic music and phase music. They are both generated by a small number of finite repeating measures, which result in complex polyrythmic sounds emerging from the combination. Quasiperiodic music and phase music both feature musical parts drifting in and out of sync. In phase music, this occurs through the gradual slowing of the tempo of one part of the music, and in quasiperiodic music this occurs as the starting point of two parts get closer and closer together after each quasiperiod. In some sense, quasiperiodic music might be viewed as an alternative approach to phasing.

## 6. Conclusion

We used a common mathematical definition of quasiperiodic function, taken from Avron and Simon (1981). We then applied the structure of a quasiperiodic function to music, describing an original music composition (Ong, 2020). Then, we contextualized this creative experiment within recent studies on mathematical music theory. Finally, we made comparisons between quasiperiodic music and related ideas in music composition and theory, for example, algorithmic music (with, in particular, Lindenmayer systems), and phase music.

## Acknowledgments

This work was partially supported by a grant from the Fundamental Research Grant Scheme from the Malaysian Ministry of Education (Grant No: FRGS/1/2018/STG06/XMU/02/1) and two Xiamen University Malaysia Research Funds (Grant Numbers: XMUMRF/2018-C1/IMAT/0001 and XMUMRF/2020-C5/IMAT/0011). The author also wishes to thank Jake and Jamie Fillman for many helpful discussions, Maria Mannone for creating that delightful design for the graphical abstract, as well as the anonymous reviewers for many comments that improved the content and presentation of the paper.

## Disclosure statement

No potential conflict of interest was reported by the author(s).

## Funding

This work was partially supported by a grant from the Fundamental Research Grant Scheme from the Malaysian Ministry of Education (Grant No: FRGS/1/2018/STG06/XMU/02/1) and two Xiamen University Malaysia Research Funds (Grant Numbers: XMUMRF/2018-C1/IMAT/0001 and XMUMRF/2020-C5/IMAT/0011).

## ORCID

*Darren C. Ong* http://orcid.org/0000-0003-4942-4376